\documentclass[a4paper,12pt]{amsart}
\usepackage{amssymb}
\usepackage{ifthen}
\usepackage{graphicx}
\usepackage{float}
\usepackage{caption}
\usepackage{subcaption}
\usepackage{cite}
\usepackage{amsfonts}
\usepackage{amscd}
\usepackage{amsxtra}
\usepackage{mathrsfs}
\usepackage[usenames]{color}

\newtheorem{thm}{Theorem}
\newtheorem{cor}{Corollary}
\newtheorem{lem}{Lemma}

\newtheorem{rem}{Remark}
\newtheorem{example}{Example}
\newtheorem{defn}{Definition}
\newtheorem{prob}{Problem}

\newtheorem{conj}{Conjecture}
\theoremstyle{definition}

\newcounter {own}
\def\theown {\thesection  .\arabic{own}}

\newenvironment{pf}[1][]{%
 \vskip 3mm
 \noindent
 \ifthenelse{\equal{#1}{}}%
  {{\slshape Proof. }}%
  {{\slshape #1.} }%
 }%
{\qed\bigskip}

\newcounter{alphabet}
\newcounter{tmp}
\newenvironment{Thm}[1][]{\refstepcounter{alphabet}%
\bigskip%
\noindent%
{\bf Theorem \Alph{alphabet}}%
\ifthenelse{\equal{#1}{}}{}{ (#1)}%
{\bf .} \itshape}{\vskip 8pt}

\makeatletter
\newcommand{\Ref}[1]{\@ifundefined{r@#1}{}{\setcounter{tmp}{\ref{#1}}\Alph{tmp}}}
\makeatother

\newcommand{\IR}{{\mathbb R}}

\newcommand{\IC}{{\mathbb C}}
\newcommand{\ID}{{\mathbb D}}

\def\be{\begin{equation}}
\def\ee{\end{equation}}

\newcommand{\bee}{\begin{enumerate}}
\newcommand{\eee}{\end{enumerate}}

\newcommand{\blem}{\begin{lem}}
\newcommand{\elem}{\end{lem}}
\newcommand{\bthm}{\begin{thm}}
\newcommand{\ethm}{\end{thm}}
\newcommand{\bcor}{\begin{cor}}
\newcommand{\ecor}{\end{cor}}
\newcommand{\beg}{\begin{example}}
\newcommand{\eeg}{\end{example}}
\newcommand{\begs}{\begin{examples}}
\newcommand{\eegs}{\end{examples}}
\newcommand{\bdefe}{\begin{defn}}
\newcommand{\edefe}{\end{defn}}
\newcommand{\bprob}{\begin{prob}}
\newcommand{\eprob}{\end{prob}}
\newcommand{\bei}{\begin{itemize}}
\newcommand{\eei}{\end{itemize}}

\newcommand{\bcon}{\begin{conj}}
\newcommand{\econ}{\end{conj}}
\newcommand{\bcons}{\begin{conjs}}
\newcommand{\econs}{\end{conjs}}
\newcommand{\bprop}{\begin{propo}}
\newcommand{\eprop}{\end{propo}}
\newcommand{\br}{\begin{rem}}
\newcommand{\er}{\end{rem}}
\newcommand{\brs}{\begin{rems}}
\newcommand{\ers}{\end{rems}}
\newcommand{\bo}{\begin{obser}}
\newcommand{\eo}{\end{obser}}
\newcommand{\bos}{\begin{obsers}}
\newcommand{\eos}{\end{obsers}}
\newcommand{\bpf}{\begin{pf}}
\newcommand{\epf}{\end{pf}}
\newcommand{\ba}{\begin{array}}
\newcommand{\ea}{\end{array}}
\newcommand{\beq}{\begin{eqnarray}}
\newcommand{\beqq}{\begin{eqnarray*}}
\newcommand{\eeq}{\end{eqnarray}}
\newcommand{\eeqq}{\end{eqnarray*}}

\newcommand{\ds}{\displaystyle}

\newcounter{minutes}
\divide\time by 60
\newcounter{hours}
\multiply\time by 60 \addtocounter{minutes}{-\time}

\begin{document}
\bibliographystyle{amsplain}
\title[]{Injectivity of sections of close-to-convex harmonic mappings with convex analytic part}

\thanks{
File:~\jobname .tex,
          printed: \number\year-\number\month-\number\day,
          \thehours.\ifnum\theminutes<10{0}\fi\theminutes}

%

%
\author{Anbareeswaran Sairam Kaliraj}
\address{A. Sairam Kaliraj, Indian Statistical Institute (ISI), Chennai Centre, Old No. 110, Chateu de Ampa, First floor, Nelson Manickam Road, Aminjikarai, Chennai 600 029, India.}
\email{sairamkaliraj@gmail.com}
%


\date{\today  
}

\begin{abstract}
In this article, we determine two point distortion theorem and sharp coefficient estimates for the families of close-to-convex harmonic mappings whose analytic part is a convex function of order $\alpha$. By making use of these results, we determine the radius of univalence of sections of these families in terms of zeros of certain equation. Lower bound for the radius of univalence has been obtained explicitly for the case $\alpha = 1/2$. Comparison of radius of univalence of the sections have been shown by providing a table of numerical estimates for the special choices of $\alpha$.
\end{abstract}
\thanks{ }

\maketitle
\pagestyle{myheadings}
\markboth{A. Sairam Kaliraj}{Injectivity of sections of  close-to-convex harmonic mappings with convex analytic part}
\section{Introduction and Motivation}\label{Sai11Sec1}
Harmonic mappings are useful in the study of fluid flow problems (see \cite{Aleman_Constantin}). Furthermore, univalent harmonic functions having special geometric 
property such as convexity, starlikeness, and close-to-convexity  arises naturally while dealing with planar fluid dynamics problems. For example, in \cite[Theorem 4.5]{Aleman_Constantin},  
Aleman et al. considered a fluid flow problem on a convex domain $\Omega_0$ satisfying an interesting geometric property. Furthermore, the harmonic mappings which appears as a solution to
some of these real world problems are very complicated and hence evaluating the values of such functions are challenging. In this connection, it is interesting to consider the problem of 
approximating harmonic mappings by harmonic polynomial without compromising the univalency of the polynomial. With this brief motivation, we shall begin considering the partial sums for the 
univalent analytic functions and univalent harmonic mappings. \\

Let $\mathcal{S}$ denote the class of all normalized univalent analytic functions $\phi$ defined in the unit disk
$\mathbb{D} = \{z: \; |z|<1\}$ such that
\begin{equation}\label{ser_rep_ana}
\phi(z)= z + \sum_{k=2}^{\infty} a_kz^k.
\end{equation}
For $n \geq 2$, the $n^{th}$ partial sum of $\phi$ is the polynomial defined by
$$
s_n(\phi)(z)= z + \sum_{k=2}^{n} a_kz^k.
$$
In \cite{Szego} Szeg\"{o} proved that the partial sum $s_n(\phi)$ is univalent in $|z| < 1/4$ for all $\phi\in\mathcal{S}$ and $n \geq 2$.
In \cite{Rober41}, Robertson proved that the $n^{th}$ partial sum of the Koebe function is
starlike in the disk $|z|<1-3n^{-1}\log n$, for $n \geq 5$, and that $3$ can not be replaced by smaller constants.  The general theorems on convolutions \cite{Rusch_Sheil} (see also \cite[P.256, 273]{Duren}) allow one to infer that $s_n(\phi)$ is convex, starlike, or close-to-convex in the disk $|z|<1-3n^{-1}\log n$,
for $n \geq 5$, whenever $\phi$ is convex, starlike, or close-to-convex in $\ID$ (for details on the class $\mathcal{S}$ and its geometric subclasses one can refer \cite{Duren}), respectively.  However, the exact radius of univalence $r_n$ of $s_n(\phi)$ remains an open problem, if $\phi \in \mathcal{S}$. Jenkins \cite{Jenkins} proved that $s_n(\phi)$ is univalent in $|z|<r_n$ for $\phi\in\mathcal{S}$, where the radius of univalence $r_n$ is atleast $1-(4\log n - \log(4\log n))/n$ for $n \geq 8$.
However, by making use of the exact coefficient bounds, one could get $r_n >  1-(4\log n - 2\log(\log n))/n$ for $n \geq 7$. Moreover, $1-(4\log n - 2\log(\log n))/n > 1-(4\log n - \log(4\log n))/n$ for $n \geq 55$.\\

Denote by ${\mathcal H}$ the class of all complex-valued harmonic functions $f=h+\overline{g}$ in
${\mathbb D}$ normalized by $h(0)=g(0)=0=h'(0)-1 $. We call $h$ and $g$, the analytic and the co-analytic
parts of $f$, respectively, and obviously they have the following power series representation
\be\label{PSSerRep}
h(z)=z+\sum _{k=2}^{\infty}a_k z^k ~\mbox{ and }~ g(z)=\sum _{k=1}^{\infty}b_k z^k, ~z \in \mathbb{D}.
\ee
Throughout the discussion we shall use this representation.
The Jacobian $J_f$ of $f=h+\overline{g}$ is $J_f(z) = |h'(z)|^2-|g'(z)|^2$. We say that $f$ is sense-preserving in $\ID$ if $J_f(z)>0$ in $\ID$. Let $\mathcal{S}_H$ denote the class of all sense-preserving
harmonic univalent mappings $f \in {\mathcal H}$ and set $\mathcal{S}^0_H=\{f \in \mathcal{S}_H:\,  f_{\overline{z}}(0)=0\}$.
For many basic results on univalent harmonic mappings in  $\mathcal{S}^0_H$ and its well-known geometric
subclasses, namely, ${\mathcal K}_H^0$, ${\mathcal S}_H^{*0}$, and ${\mathcal C}_H^0$
mapping $\ID$ onto, convex, starlike, and close-to-convex domains, respectively, we refer to \cite{Duren:Harmonic,PonRasi2013}. Here we recall a new subclass of $\mathcal{S}^0_{H}$ namely, $\mathcal{S}^0_H(\mathcal{S})$ introduced in \cite{SaiPIAS}, where
$$\mathcal{S}^0_H(\mathcal{S}) = \left\{h+\overline{g} \in \mathcal{S}^0_H :\,
h+e^{i \theta}g \in \mathcal{S}~ \mbox{for some}~~\theta \in \mathbb{R} \right\}.$$

For $f=h+\overline{g} \in \mathcal{S}^0_H$ with power series representation as in \eqref{PSSerRep}, the sections/partial sums
$s_{n,m}(f)$ of $f$ are defined as
$$s_{n,m}(f)(z)=s_n(h)(z)+\overline{s_m(g)(z)},
$$
where $n\geq 1$ and $m\geq 2$. However, the special case $m=n\geq 2$ seems interesting in its own merit.
From the above definition, it is clear that partial sums of $f$ can be thought of as an approximation of $f$ by complex-valued harmonic polynomials and thus,
approximation of univalent harmonic mappings by univalent harmonic polynomials might lead to new applications. For fundamental results on the partial sums of univalent harmonic mappings,
one can refer to \cite{LiSamyNA1, LiSamyNA2, PonSai8, PonSai9}. We recall a few results from \cite{PonSai8},
which are motivation for the problem which we consider in this article.

\begin{Thm}\label{RPS7Thm2}{\rm \cite[Theorem 1]{PonSai8}}
Let $f=h+\overline{g} \in \mathcal{S}^0_H$ with the series representation as in \eqref{PSSerRep}. Suppose that $f$ belongs to ${\mathcal C}_H^0$, the class of close-to-convex harmonic mappings or $\mathcal{S}^0_H(\mathcal{S})$. Then the section $s_{n,m}(f)$ is univalent in the disk $|z|<r_{n,m}$. Here $r_{n,m}$ is the unique positive
root of the equation $\psi(n,m,r)=0$, where
$$
\psi(n,m,r) = \frac {1}{12 r}
\left(\frac {1-r}{1+r}\right)^3 \left[1-\left(\frac {1-r}{1+r}\right)^{6}\right] - R_n - T_m,
$$
with
$$
R_n = \sum_{k=n+1}^{\infty}A_k r^{k-1}, ~~ T_m = -\sum_{k=m+1}^{\infty}A_{-k} r^{k-1}, ~\mbox{ where }~
A_k=\frac{k(k+1)(2k+1)}{6},
$$

In particular, every section $s_{n,n}(f)(z)$ is univalent in the disk $|z|<r_{n,n}$, where
$$r_{n,n} > r^L_{n,n}:=1- \frac{(7\log n  - 4\log(\log n))}{n} ~\mbox{ for }~ n \geq 15.
$$
Moreover, $r_{n,m} \geq r^L_{l,l}$,  where  $l=\min\{n,m\}\geq 15$.
\end{Thm}

For functions in the convex family $\mathcal{K}^0_H$ of harmonic mappings, we have the following interesting result, in which the lower bounds on $r_{n,m}$ is better than that of the bounds in Theorem \Ref{RPS7Thm2}.

\begin{Thm}\label{RPS7Thm3}{\rm \cite[Theorem 2]{PonSai8}}
Let $f=h+\overline{g} \in \mathcal{K}^0_H$ with the series representation as in \eqref{PSSerRep}. Then the section $s_{n,m}(f)$ is univalent in the disk $|z|<r_{n,m}$, where $r_{n,m}$ is the unique positive root of the equation $\mu(n,m,r)=0$. Here
\be\label{PS7_eq5}
\mu(n,m,r)=\frac {1-r}{(1+r)^3} - \sum_{k=n+1}^{\infty}\left[ \frac{k(k+1)}{2} r^{k-1} \right] - \sum_{k=m+1}^{\infty}\left[ \frac{k(k-1)}{2} r^{k-1} \right].
\ee
In particular, for $n\geq5$, and $\theta \in \IR$, the harmonic function
$$s_{n,n}(f;\theta)(z)=s_n(h)(z)+e^{i\theta}\overline{s_n(g)(z)}
$$
is univalent and close-to-convex in the disk $|z| < 1-3n^{-1}\log n $.
Moreover, we have $r_{n,m} \geq 1-(4\log l - 2\log(\log l))/l$, where $l=\min\{n,m\} \geq 7 $.
\end{Thm}

For $\alpha \geq -1/2$, set 
$${\mathcal F}(\alpha)=\left\{f=h+\overline{g}\in {\mathcal H} :\, {\rm Re\,} \left (1 + \frac{zh''(z)}{h'(z)}\right ) >
 \alpha, g'(z)=e^{i\theta}zh'(z) \mbox{ for all }~z\in\ID \right\}.
$$

\br Let us consider the family of close-to-convex harmonic functions defined by ${\mathcal F} = {\mathcal F}(-1/2)$.
A keen observation of the proof of Theorem \Ref{RPS7Thm3} reveals that the results in Theorem \Ref{RPS7Thm3} is valid for the class
${\mathcal F}$ too (It should be noted that the members of ${\mathcal F}$ are not necessarily convex). Comparison of this fact with the result in Theorem \ref{RPS7Thm2}, motivates us to consider classes of close-to-convex harmonic mappings with special condition on the analytic part of $f$.
\er


In this article, we determine the radius of univalence $r_{n,m}$ of $s_{n,m}(f)$, when $f \in {\mathcal F}(\alpha)$. Lower bound for $r_{n,m}$ are determined for the special choice of $\alpha$. The paper is organized as follows. In Section \ref{Sai11Sec2}, we recall a few known results from the literature, which are helpful in the proof of our main theorems. In Section \ref{Sai11Sec3}, we present our main Theorems.

\section{Useful Results}\label{Sai11Sec2}

The following result due to Bazilevich \cite{Bazilevich} gives the necessary and sufficient condition for a normalized analytic function to be univalent in $\ID$.

\begin{Thm}{\rm \cite{Bazilevich}}\label{uni_nec_suf_Analytic}
An analytic function $\phi$ defined in $\ID$ and determined by \eqref{ser_rep_ana} is univalent in $\ID$ if and only if for each
$z \in \ID$ and each $t\in[0, \pi/2]$,
\be\label{PS7inteq1}
\frac {\phi(re^{i\eta})-\phi(re^{i\psi})} {re^{i\eta}-re^{i\psi}} :=\sum_{k=1}^{\infty}a_k \frac{\sin kt}{\sin t} z^{k-1} \ne 0,
\ee
where $t=(\eta-\psi)/2$, $z=re^{i(\eta+\psi)/2}$ and $\left .\frac{\sin kt}{\sin t}\right |_{t=0}=k$.
\end{Thm}

In \cite{Starkov},  Starkov  generalized this result to the class of normalized sense-preserving harmonic mappings in the following form.

\begin{Thm}{\rm \cite{Starkov}}\label{uni_nec_suf}
A sense-preserving harmonic function $f=h+\overline{g}$ defined in $\ID$ determined by \eqref{PSSerRep} is univalent in $\ID$ if and only if for each
$z \in \ID\setminus\{0\}$ and each $t\in(0, \pi/2]$,
\be\label{PS7inteq2}
 \frac {f(re^{i\eta})-f(re^{i\psi})} {re^{i\eta}-re^{i\psi}}  :=\sum_{k=1}^{\infty}\left[ (a_k z^k - \overline{b_k z^k})\frac{\sin kt}{\sin t} \right] \ne 0,
\ee
where $t=(\eta-\psi)/2$ and $z=re^{i(\eta+\psi)/2}$.
\end{Thm}

The following two point distortion theorem plays a crucial role in the proof of our main results.

\begin{Thm}\label{two_po_dis}
If $f \in \mathcal{S}$, $r \in (0, 1)$, $t, \psi \in \IR$, then
$$\left|\frac {f(re^{it})-f(re^{i\psi})} {re^{it}-re^{i\psi}}\right|\ge \frac {1-r^2}{r^2}
|f(re^{it})| \, |f(re^{i\psi})|.
$$
\end{Thm}

\section{Distortion Theorem and Partial sums problem}\label{Sai11Sec3}

By making use of Theorem \Ref{two_po_dis}, we derive two point distortion theorem for functions in the family
${\mathcal F}(\alpha)$. This result is very crucial in the proof of our main theorem.

\bthm\label{Sai11thm1}{\bf (Two point distortion Theorem)}
Suppose that $f=h+\overline{g} \in {\mathcal F}(\alpha)$ for some  $\alpha$~ $(0 \leq \alpha < 1)$. For 
each $\lambda \in \IC$ such that $|\lambda| =1$ define 
$$F_{\lambda}(z):= h(z) + \lambda g(z)
$$
Then, for any $t, \psi \in \IR$ such that $t \neq \psi$, $f$ 
satisfies the following inequality
\beq\label{Har_dist}
\left|\frac {F_{\lambda}(re^{it})-F_{\lambda}(re^{i\psi})} {re^{it}-re^{i\psi}}\right| &\ge& \frac{1-r^2}{r^2} L^2(r,\alpha) =: A(r, \alpha)~~\mbox{ for }~ 0 < r < 1,\\
\left|\frac {f(re^{it})-f(re^{i\psi})} {re^{it}-re^{i\psi}}\right| &\ge&  A(r, \alpha)~~\mbox{ for }~ 0 < r < 1, \nonumber
\eeq
where
$$L(r,\alpha) := \left\{
    \begin{array}{ll}
        \ds \frac{2r}{1+r}-\log(1+r)  & \mbox{if } \alpha=0 \\
        \ds 2\log(1+r)-r & \mbox{if } \alpha=1/2\\
        \ds \frac{(1+r)^{2\alpha}(1+r+2\alpha-2r\alpha)-(1+r)(1+2\alpha)}{2\alpha(1+r)(2\alpha -1)} & \mbox{if }
        \alpha \in (0, 1) \setminus \{ 1/2\}.
    \end{array} \right.
$$
\ethm

\bpf
Let $f=h+\overline{g} \in {\mathcal F}(\alpha)$ for some  $\alpha$ $(0 \leq \alpha < 1)$ and $F_{\lambda}(z):= h(z) + \lambda g(z)$.
From the definition of $f \in {\mathcal F}(\alpha)$, it is clear that $h$ is a convex function of order $\alpha$ and
$g'(z) = e^{i\theta} z h'(z)$. Therefore, from a result of Clunie and Sheil-Small \cite[Theorem 5.17]{Clunie-Small-84}, it is
clear that $F_{\lambda}$ is close-to-convex in $\ID$ for all $\lambda $ such that $|\lambda| =1$. For every pair of points $re^{it}$ and $re^{i\psi}$, we can find a $\lambda$ such that
$$(h(re^{it})-h(re^{i\psi}))+\overline{(g(re^{it})-g(re^{i\psi}))} = (h(re^{it})-h(re^{i\psi}))+\lambda(g(re^{it})-g(re^{i\psi}))
$$
Therefore
\beq\label{Sai11_th1eq1}
\left|\frac {f(re^{it})-f(re^{i\psi})} {re^{it}-re^{i\psi}}\right| &=& \left|\frac {F_{\lambda}(re^{it})-F_{\lambda}(re^{i\psi})} {re^{it}-re^{i\psi}}\right| \nonumber \\
&\ge& \frac {1-r^2}{r^2}
|F_{\lambda}(re^{it})|\, |F_{\lambda}(re^{i\psi})|. ~~~\mbox{(by Theorem \Ref{two_po_dis})}
\eeq
In order to complete the proof, we need to find the lower bounds of $|F_{\lambda}(z)|$.

Let $\gamma$ be the preimage under $F_{\lambda}$ of the line segment joining $0$ and $F_{\lambda}(z)$. Then
\beq\label{Sai11_th1eq2}
|F_{\lambda}(z)|&=& \int\limits_\gamma \left | F'_{\lambda}(\zeta)\right| |d\zeta|\,\nonumber\\
& \geq & \int\limits_\gamma (|h'(\zeta)|-|g'(\zeta)|) |d\zeta|\nonumber\\
& \geq & \int\limits_\gamma (1-|\zeta|)|h'(\zeta)|\, |d\zeta| \nonumber\\
& \geq & \int\limits_0^r \frac{(1-\rho)}{(1+\rho)^{2(1-\alpha)}}\, d\rho,
\eeq
where, the last inequality is the consequence of the following well known distortion
theorem for the convex function of order $\alpha ~(0 \leq \alpha < 1)$ (see \cite[p.~139, Vol. 1, Theorem 1]{Go})
$$
|h'(z)| \ge \frac{1}{(1+r)^{2(1-\alpha)}}, ~ r = |z| < 1.
$$
The desired conclusion follows if we use the inequality in \eqref{Sai11_th1eq2} in \eqref{Sai11_th1eq1}.
\epf

Next, we provide the sharp coefficient estimates for the functions in the family ${\mathcal F}(\alpha)$. 

\bthm\label{Sai11thm2}
Suppose that $f=h+\overline{g} \in {\mathcal F}(\alpha)$ for some  $\alpha$~ $(0 \leq \alpha < 1)$ with the series representation as in \eqref{PSSerRep}. Then, for all $n \ge 2$, the coefficients of $f$ satisfy the following inequality
\be\label{Sai11_th2eq3}
|a_n| \le A_n(\alpha) ~\mbox{ and }~ |b_n| \le \frac{n-1}{n-2\alpha}A_n(\alpha), 
\ee
where 
\be\label{Sai11_th2eq3.1}
A_n(\alpha) = \frac{1}{n!}\prod_{j=2}^n (j-2\alpha).
\ee
All these bounds are sharp and the equality in each inequality is attained for the close-to-convex functions
$f_{\alpha}(z)=h_{\alpha}(z)+\overline{g_{\alpha}(z)}$ and its rotations, where
\be\label{Sai11_th2eq4}
f_{\alpha}(z)=
 \left \{
 \begin{array}{ll}
  \ds \frac{1-(1-z)^{2\alpha -1}}{2\alpha - 1}+\overline{\frac{1-(1-z)^{2\alpha -1}(1+z(2\alpha -1))}
  {2\alpha(2\alpha - 1)}}  & \mbox{if $\alpha\neq0, 1/2$},\\ \vspace{0.2cm}
  \ds \frac{z}{1-z} + \overline{\frac{z}{1-z}+\log(1-z)} & \mbox{if $\alpha=0$},\\
  \vspace{0.2cm}
  \ds -\log(1-z)-\overline{(z+\log(1-z))} & \mbox{if $\alpha=1/2$}.
  \end{array}
 \right.
\ee
\ethm
\bpf
The proof follows from the coefficient bounds for the convex functions of order $\alpha$ (see \cite[p.~140, Vol. 1, Theorem 2]{Go}).
\epf

\begin{thm}\label{Sai11thm3}
Let $f=h+\overline{g} \in {\mathcal F}(\alpha)$ with the series representation as in \eqref{PSSerRep}. Then for $\theta \in \IR$, the harmonic function
$$s_{n,m}(f;\theta)(z)=s_n(h)(z)+e^{i\theta}\overline{s_m(g)(z)}
$$
univalent in the disk $|z|<r_{n,m}$, where $r_{n,m}$ is the unique positive root of the equation $\mu(n,m,r, \alpha)=0$ in $(0, 1)$. Here
\be\label{Sai11_th2eq5}
\psi(n,m,r,\alpha)=A(r, \alpha) - \sum_{k=n+1}^{\infty}\left[ k A_k(\alpha) r^{k-1} \right] - \sum_{k=m+1}^{\infty}\left[ \frac{k(k-1)}{(k-2\alpha)}A_k(\alpha) r^{k-1} \right],
\ee
where $A(r, \alpha)$ and $A_k(\alpha)$ are defined as in \eqref{Har_dist} and \eqref{Sai11_th2eq3}, respectively. When $\alpha =1/2$,
we have $r_{n,m} \geq 1-(2\log l)/l$, where $l=\min\{n,m\} \geq 3$.
\end{thm}

\bpf
Suppose that $f=h+\overline{g} \in {\mathcal F}(\alpha)$ with the series representation as in \eqref{PSSerRep}. For $\theta \in \IR$, we set
$$F_{r, \theta}(z)=\frac{h(rz)}{z} + e^{i\theta} \frac{g(rz)}{z}$$ for $0 < r < 1$, one could see that
$$ F_{r,\theta}(z) = z+\sum _{k=2}^{\infty}a_k r^{k-1} z^k + e^{i\theta} \sum _{k=2}^{\infty}{b_k r^{k-1} z^k}.
$$
We shall prove that $s_{n,m}(F_{r, \theta})(z)$ is univalent in $\ID$ for all values of $\theta$. This would imply that $s_{n,m}(f;\theta)$ is univalent in the disk $|z|< r$
for all values of $\theta$. By making use of  Theorem \Ref{uni_nec_suf_Analytic}, we see that $s_{n,m}(F_{r, \theta})(z)$
is univalent in $\ID$ if and only if the associated section $P_{n,m,r,\alpha}(z)$
has the property that
$$P_{n,m,r,\alpha}(z) := \sum_{k=1}^{M}\left[ (a'_k z^k + {b'_k z^k})\frac{\sin kt}{\sin t} \right] \ne 0, ~\mbox{ for all }~ z\in\ID\setminus\{0\}, ~\mbox{ and }~
t\in[0,\pi/2],
$$
where $M = \max\{n, m\}$, $l=\min\{n, m\}$, $a'_k = a_k r^{k-1}$, $b'_k = b_k r^{k-1}$ for all $k \leq l$,
$$
a'_k= \left \{ \begin{array}{lr}
\ds a_k r^{k-1} &  \mbox{for all $k > l$ if $M = n$},\\
0 & \mbox{for all $k > l$ if $M > n$},
\end{array}
\right.
$$
and
$$
b'_k=\left \{ \begin{array}{lr}
\ds e^{i\theta} b_k r^{k-1} &  \mbox{for all $k > l$ if $M = m$},\\
0 & \mbox{for all $k > l$ if $M > m$}.
\end{array}
\right.
$$
Setting $t=(\eta-\psi)/2$, $z=\rho e^{i(\eta+\psi)/2}\in\ID$ in \eqref{PS7inteq2} and from the univalency of $F_{r, \theta}$ for $0<r < 1$, we get that
\be\label{Sai11_th2eq6}
\left|\sum_{k=1}^{\infty}\left[ (a_k z^k +  e^{i\theta}{b_k z^k})r^{k-1}\frac{\sin kt}{\sin t} \right]\right| \geq  A(r, \alpha).
\ee
In order to find a lower bound for $|P_{n,m,r,\alpha}(z)|$, we need to find an upper bound for
$$\left|R_{n,m,r,\alpha}(z)\right| = \left| \sum_{k=n+1}^{\infty}\left[ a_k r^{k-1} z^k \frac{\sin kt}{\sin t} \right] + \sum_{k=m+1}^{\infty}\left[ {e^{i\theta}(b_k r^{k-1} z^k)} \frac{\sin kt}{\sin t} \right]\right|.
$$

From Theorem \ref{Sai11thm2}, it follows that
\beq\label{Sai11_th2eq7}\nonumber
|R_{n,m,r,\alpha}(z)| &\leq& \sum_{k=n+1}^{\infty}\left[ k A_k(\alpha) r^{k-1} \right] + \sum_{k=m+1}^{\infty}\left[ \frac{k(k-1)}{(k-2\alpha)}A_k(\alpha) r^{k-1} \right], \\
 &=& R_n (\alpha) + T_m(\alpha),
\eeq
where $A(r, \alpha)$ and $A_k(\alpha)$ are defined as in \eqref{Har_dist} and \eqref{Sai11_th2eq3}, respectively.
From \eqref{Sai11_th2eq6} and \eqref{Sai11_th2eq7}, we get that
$$|P_{n,m,r,\alpha}(z)| \geq A(r, \alpha) - R_n (\alpha) - T_m (\alpha)=\psi(n,m,r, \alpha).
$$
The inequality $P_{n,m,r,\alpha}(z) \neq 0$ holds for all $z\in\ID\setminus\{0\}$,  whenever $\psi(n,m,r,\alpha)>0$, where $\psi(n,m,r,\alpha)$ is
defined by \eqref{Sai11_th2eq5}. This gives that  $\psi(n,m,r,\alpha)>0$ for all $r \in (0, r_{n,m})$, where $r_{n,m}$ is the positive root of the equation $\psi(n,m,r,\alpha)=0$, which lies in the interval $(0, 1)$. This observation proves that $ F_{r, \theta}$ is univalent in $\ID$ for all $\theta \in \IR$, which implies that $s_{n,m}(f;\theta)$ is univalent (see \cite{Rodri-Maria}) in the disk $|z|< r_{n,m}$ for all $\theta \in \IR$. This completes the proof of the first part of the theorem.

From the above discussion, it is apparent that $r_{n,m} \geq r_{l,l}$, where $l=\min\{n, m\}\geq 2$.
Next, we shall consider the special case $m=n$ and determine the lower bound for $r_{n,n}$ with certain restriction on $n$ and special choices of
$\alpha$.

Let us consider the case $\alpha = 1/2$. In this case, the sufficient condition \eqref{Sai11_th2eq5} for the univalence of $s_{n,n}(f;\theta)$ reduces to
\beqq\label{PS7_eq4}
\psi(n,n,r,1/2)&=&A(r, 1/2) - R_n (1/2) - T_n (1/2)\\
&=& (1-r^2) \left(\frac{2\log(1+r)-r}{r}\right)^2 - 2\sum_{k=n+1}^{\infty} r^{k-1} \\
&\ge& \frac{1-r}{(1+r)^2} - \frac{2r^n}{(1-r)}
\eeqq
From the fact that $\lim_{n\rightarrow\infty}\psi(n,n,r,1/2) > 0 $, for all $r$ such that $0 < r< 1$,
it is clear that the radius of univalence $r_{n,n}\rightarrow 1$ as $n\rightarrow \infty$.
Setting $r = 1-x/n$, where $x=o(n)$, we see that
$$
\psi(n,n,r,1/2) = \frac{(1-r)^2 - 2(1+r)^2r^n}{(1-r)(1+r)^2} > 0
$$
whenever $0<t(x, n)<1$, where
$$t(x, n)=2e^{-x} \left(\frac{2n}{x} - 1\right)^2.
$$
Now, we shall prove that $r_{n,n}>1- \gamma_n/n$ and we shall explicitly determine the value $\gamma_n$ for large values of $n$. In order to prove that, it is sufficient to prove that $t(x, n)$ is a decreasing function in $x$, whenever
$$\gamma_n \leq  x \leq n,~~0 < t(\gamma_n, n) < 1 ~\mbox{ and }~ t(n, n)>0.
$$
Since $e^{x}t(x, n)= O(n^2)$, we may set $\gamma_n=2\log n$. It is easy to see that $1- \gamma_n/n > 0$ and increasing for all $n \geq 3$.
For $n \geq 3$, we shall prove that $r_{n,n}>1- \gamma_n/n$.

For every fixed integer  $n \ge 3$, let us consider $x \in [\gamma_n, n]$ and we prove that  $t(x,n)$ decreasing function with respect to $x$.
From the definition of $t(x,n)$, it is clear that it is a product of two positive decreasing functions. Hence $t(x,n)$ is a decreasing function with respect to $x$ for every fixed integer $n \ge 3$. Further, it is easy to see that $t(n, n) < 1$ for all $n \ge 3$. In order to complete the proof, it is enough to show that $t(\gamma_n, n) < 1$ for all $n \ge 3$.

\beqq
t(\gamma_n, n) &=& \frac{2}{n^2}\left(\frac{n}{\log n} - 1 \right)^2\\
&=& 2 \left(\frac{1}{\log n} - \frac{1}{n} \right)^2\\
&<& 1 ~~\mbox{ for all }~~ n \ge 3.
\eeqq

The proof is complete.
\epf

\br
From the above computation, it is clear that one could estimate the lower bound for $r_{n,m}$, whenever $\alpha \in [0, 1)$. In order to compare the results on the radius of univalence of the 
partial sums, when we consider functions $f \in {\mathcal F}(\alpha)$ for various values of $\alpha$, we provide the following Corollary.
\er

\bcor
Suppose that $f \in {\mathcal F}(\alpha)$. Then the value of $n$ for which $s_{n,n}(f;\theta)(z)$ is univalent in the disk $|z| < \rho$ is
formulated in Table \ref{tab2}:
\ecor
\begin{table}[H]
\begin{center}
\begin{tabular}{|l|l|l|l|l|l|l|}
  \hline
 {\bf Value of $\rho$}   & ${\bf \alpha =-1/2}$ & ${\bf \alpha =0}$  & ${\bf \alpha =1/4}$ & ${\bf \alpha =1/2}$  & ${\bf \alpha =3/4}$ \\
  \hline
  {\bf 1/4 }              &   $n \ge 4$     &   $n \ge 3$    &   $n \ge 2$   &  $n \ge 2$   &  $n \ge 2$     \\
  \hline
  {\bf 1/2 }             &   $n \ge 12$    &   $n \ge 8$    &  $n \ge 6$    &  $n \ge 4$   &  $n \ge 3$      \\
  \hline
  {\bf 3/4 }             &   $n \ge 43$    &   $n \ge 29$   &  $n \ge 22$   &  $n \ge 16$  &  $n \ge 9$       \\
  \hline
  {\bf 9/10}             &   $n \ge 160$   &   $n \ge 111$  &  $n \ge 86$   &  $n \ge 61$  &  $n \ge 35$       \\
  \hline
\end{tabular}
\end{center}
\caption{Values of $n$ for which $s_{n,n}(f;\theta)(z)$ is univalent in the disk $|z| < \rho$ \label{tab2}}
\end{table}
\bpf
The above lower bounds have been estimated with the help of Mathematica software by estimating $\psi(n,n,r,\alpha)=0$ in the equation \eqref{Sai11_th2eq5}, corresponding to the values of $r=\rho$, fixing $\alpha=0, 1/4, 1/2, 3/4$. The estimates for the case 
$\alpha = -1/2$ have been obtained by estimating $\mu(n,n,r) =0 $ in \eqref{PS7_eq5} by fixing $r = \rho$.
\epf

\subsection*{Acknowledgements} The author thanks National Board for Higher Mathematics (which comes under the Department of Atomic Energy, Government of India) for supporting this research work by providing NBHM postdoctoral fellowship.

\end{document}